\documentclass[conference]{IEEEtran}

\usepackage{ifpdf}
\usepackage{cite}

\usepackage{xcolor}
\ifCLASSINFOpdf
  \usepackage[pdftex]{graphicx}
\fi
\usepackage{amsmath}
\usepackage{algorithmic}
\usepackage{bm}
\usepackage{hyperref}

\newtheorem{Definition}{Definition}

\begin{document}

\title{Congestion management within a multi-service scheduling coordination scheme for large battery storage systems}



\author{\IEEEauthorblockN{Cl{\'e}mentine Straub \IEEEauthorrefmark{1}\IEEEauthorrefmark{2},
Jean Maeght\IEEEauthorrefmark{2},
Camille Pache\IEEEauthorrefmark{2},
Patrick Panciatici\IEEEauthorrefmark{2} and
Ram Rajagopal\IEEEauthorrefmark{3}}
\IEEEauthorblockA{\IEEEauthorrefmark{1}Laboratory of Signals and Systems (L2S), 
CentraleSup{\'e}lec, Universit{\'e} Paris-Saclay, Gif sur Yvette, France}
\IEEEauthorblockA{\IEEEauthorrefmark{2}French Transmission System Operator, R{\'e}seau de Transport d'{\'E}lectricit{\'e} (RTE), Paris La D{\'e}fense, France}
\IEEEauthorblockA{\IEEEauthorrefmark{3}Civil \& Environmental Engineering, Stanford University, Stanford, CA 94305, USA}
}

\maketitle

\begin{abstract}
There is a growing interest in the use of large-scale battery storage systems for grid services. This technology has been deployed in several countries to increase transmission systems capabilities and reliability. In particular, large-scale battery storage systems could be used for congestion management and the French Transmission System Operator (RTE) is currently installing 3 large batteries for 2020 at the sub-transmission grid level for this purpose. The battery operation for congestion management does not require the full storage capacities  at all times. Thus, the residual capacities can be offered to other services to increase batteries profitability. This paper presents the framework which will be used by RTE for battery operation scheduling to combine congestion management with other services by computing day-ahead bandwidths defining available storage capacities. The bandwidths represent safe domains for grid operation scheduling: as long as the battery operation is performed within these bandwidths, there will be no grid congestion or grid congestions will be managed.
\end{abstract}

\IEEEpeerreviewmaketitle

\section{Introduction}
A dramatic increase in renewable generation is foreseen in the years to come. In France, most of the deployment of wind and solar farms is expected at the sub-transmission level which may increase grid congestions locally. Indeed, the sub-transmission grid was initially sized to supply the local consumption and not to convey power generated at this level to the main transmission grid. Consequently, grid congestions may be created by the outbound flows correlated with strong renewable generation in some zones of the sub-transmission grid.
 Adding a battery storage in these zones may help congestion management: outbound flows can be reduced by charging the battery. When the congestion is over, the battery can discharge and send the stored renewable energy back to the main grid. Several papers present algorithms for congestion management using large batteries. For instance, \cite{almassalkhi2015model} presents a Model Predictive Control using large storage systems and redispatch to mitigate the effects of severe line-overloads with constraints on the line temperatures. In \cite{pahalawaththa2018}, the Australian Transmission System Operator (TSO)  identifies situations where batteries can be used to increase transmission capabilities. In \cite{hemmati2017stochastic}, an optimal energy storage scheduling for congestion management with a focus on uncertainties is proposed.

RTE, the French TSO, in charge of the transmission and sub-transmission grid, will install three large batteries by 2020 in three congested zones. A local real-time controller will use the battery in combination with renewable generation curtailment to manage grid congestions. The real-time controller is presented in \cite{8511518}. However, congestions are not so frequent and during periods where there is no congestion, the battery can be available for other services as long as it respects certain operation bounds.

\subsection{Scheduling for multi-service}
Despite the decreasing costs of storage \cite{schmidt2017future}, batteries remain very expensive and their use for congestion management alone does not seem profitable for the system. Stacking different services offered by the batteries, however, may increase the revenue and compensate their high costs. Multi-service for batteries is widely present in the literature and several papers propose approaches to combine services. For instance, \cite{wu2015energy} presents a framework for battery multi-service portfolio including energy arbitrage, distribution grid deferral and outage mitigation. \cite{namor2018control} summarizes the combination of services present in the literature and classifies them in three categories: energy arbitrage, provision of ancillary services and the achievement of local control objectives, in which congestion management is included. \cite{namor2018control} also proposes a model for battery sharing validated by real-time experiments in which each service has an allocated budget in power and energy.

This paper investigates battery scheduling for multi-service, and in particular the coordination of congestion management with other services. The framework differs from the existing literature by the definition of residual capacities not used for congestion management,  which can be offered to other services. These residual capacities are called bandwidths in the remaining of the paper.

\subsection{Definition of bandwidths}
The bandwidths define safe operational domains for batteries within which all congestions can be managed and which ensure that no additional congestion is created. More precisely, the bandwidths need to present the following information:
\begin{itemize}
	\item \textit{The minimum charge or discharge needed to solve congestion:} During a congestion, the battery can either be asked to charge (respectively discharge) partially or at full power rate. In the first case, the residual capacities can be granted to another service, whereas in the second case, the battery is not available for other purposes.
	\item \textit{An anticipation of future congestions:}  During a congestion, the battery charges (respectively discharges) power. The scheduling algorithm objective is to anticipate the congestion and empty (respectively fill) the battery before the congestion starts.
	\item \textit{Available capacities in absence of congestion:} Even in the absence of a congestion, it is essential to check that the battery operation will not create any congestion, and if it does, the bandwidths need to be reduced.
\end{itemize}

In \cite{wen2015enhanced} a similar notion is introduced: an Enhanced Security-Constrained Optimal Power Flow (OPF) in which batteries are used for post-contingencies corrective actions is presented. The minimum amount of energy a battery must be able to store or deliver for corrective actions in order to cope with all the contingencies is computed. This minimum energy is called margin and plays a similar role than the bandwidths. However, only margins in energy are considered ; margins in power are not necessary. Indeed, the objective of \cite{wen2015enhanced} is to find an operating point, which is adjusted in the case of a contingency. In the same vein, only one time step is represented as they are not dealing with scheduling and multi-service is not considered. In our scheduling process, we take the following definition for bandwidths.
\begin{Definition}
Bandwidths for battery operation are series of forecasted intervals in energy and in power such that a battery used in these intervals will be able to manage congestions and will not create any.
\label{def band}
\end{Definition}

\section{Modeling}

\subsection{Thermal ratings use in RTE}
\begin{figure}[!h]
	\centering
	\includegraphics[scale=0.55]{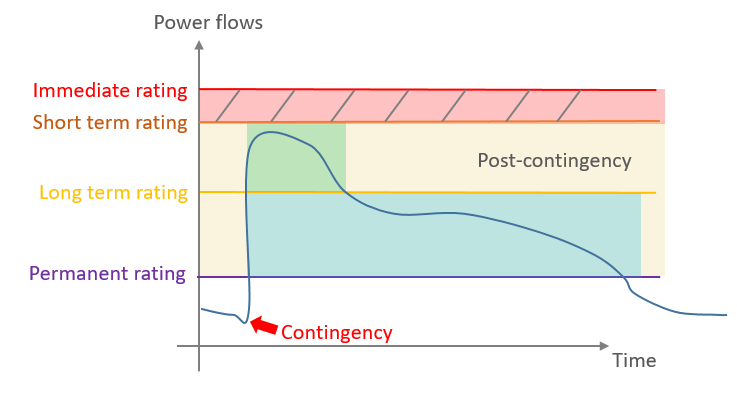}
	\caption{Thermal ratings as they are used in RTE}
	\label{rating}
\end{figure}

Bandwidths for battery operation must guarantee that transmission lines limits will not be exceeded. Several security thresholds, describing limitations under normal conditions and acceptable overloads in post-contingency states are considered in this paper. The underlying assumption for these different thresholds is that overloads will not damage equipment within a certain time frame. The thresholds allow time for automatic or human corrective actions to reschedule control after a contingency \cite{capitanescu2008new} -\cite{platbrood2011development}. The security thresholds are usually called thermal ratings. Batteries
can be used as fast corrective actions and enable to operate the grid closer to these ratings. In the scope of the paper, we consider the four thermal ratings below (see Fig. \ref{rating}):
\begin{itemize}
	\item \textit{Permanent rating:} limit under which the flow can stay indefinitely without damaging the line by overheating it.
	\item \textit{Immediate rating:} if the flow exceeds this limit, safety measures are implemented and the line is open immediately.
	\item \textit{Short term rating:} limit allowing time only for fast curative actions (batteries and topological actions).
	\item \textit{Long-term rating:} limit after fast curative actions and before long curative actions (batteries, topological changes and renewable generation curtailment).
\end{itemize}

\subsection{Bandwidths and congestions}

The objective of the bandwidths definition is to ensure that the flows will stay within the thermal ratings, both under normal conditions and in all contingency cases during operation. The bandwidths form the day-ahead battery scheduling. More precisely, they must guarantee the four following points:

\begin{itemize}
	\item \textit{Under normal conditions} the flows must not exceed the permanent rating.
	\item \textit{In all contingency cases} the flows must not exceed the immediate rating.
	\item \textit{In all contingency cases} the flows must not exceed the long-term rating after the use of fast curative actions.
	\item \textit{In all contingency cases} the flows must not exceed the permanent rating after the use of all curative actions.
\end{itemize}

Bandwidths in power and in energy are necessary to be able to guarantee these properties. The following two sections explain in details how to compute the bandwidths. First, the bandwidths in power are computed. They are then used to determine the bandwidths in energy.

\subsection{Computation of the bandwidths in power}

Bandwidths in power consist of an upper and lower bound on the admissible battery charges and discharges tolerated by the sub-transmission grid for each time step. The battery operation within these bounds ensures that ratings will not be exceeded under normal conditions and in each contingency case. Their computations rely on hourly grid forecasted states (including load, generation, topology and flows), which are already available at RTE, thanks to general day-ahead and real-time operational planning processes.

The system control variables are the battery power injection and the renewable generation curtailment. The battery power injection is decomposed in two parts, a preventive part $B$ and a curative part $B^{cur}_n$ with $n \in\mathcal{Z}^C$ the set of contingencies. The curative battery injection can indeed be adjusted depending on the contingency $n$. The renewable generation curtailment is also decomposed in a preventive part and a curative part. The preventive generation curtailment at bus $i \in \mathcal{Z}^B$, with $\mathcal{Z}^B$ the set of buses, is noted $C_i$, and the curative part is noted $C_{i,n}^{cur}$, with $n \in \mathcal{Z}^C$.

Two problems are solved for the computation of bandwidths in power: one for the determination of the lower bound, the second one for the upper bound. The two problems are very similar. They only differ by a sign in the objective function:

\begin{equation}
\label{objective}
{\text{min }} \pm B + c \cdot \sum_{i \in \mathcal{Z}^B} C_i
\end{equation}

The sign $+$ is used for the computation of the lower bound, and the sign $-$ for the upper bound. The problem can be solved independently for each time step. Consequently, indices for time steps are not present in the problem formulation. Time coupling is only considered during the computation of bandwidths in energy. The cost associated with the preventive curtailment $c$ is very high in order to prioritize the use of the battery before curtailment. An equivalent method consists in minimizing the preventive curtailment first, then fix it and minimize (and maximize) the battery injection. The curative actions terms (curative battery charge and curative curtailment) can be added in the objective function, with small weights in order to analyse the occurences of corrective actions in the results. The objective function becomes:

\begin{equation}
\label{objective2}
{\text{min }} \pm B + c_1 \cdot \sum_{i \in \mathcal{Z}^B} C_i + c_2 \cdot \sum_{n \in \mathcal{Z}^C}|B^{cur}_n| + c_3 \cdot \sum_{\substack{i \in \mathcal{Z}^B \\ n \in \mathcal{Z}^C}} C^{cur}_{i,n}
\end{equation}

In the constraints, different states of the system are considered: a state for normal conditions, as well as a state for each contingency that can occur.  Only a zone of the sub-transmission grid is represented in the model. The remaining grid is modeled with the help of Power Transfer Distribution Factors (PTDF). A PTDF is a factor which is defined for a line and a bus and  gives the amount of power going through the line if the power injection at the bus increases by 1, with a decrease of 1 in the (remote) slack bus of power flows equations (for more details, see \cite{zimmerman2011matpower} and \cite{zimmerman2010matpower}). This allows us to solve a DC OPF in the zone only and consider the flows coming from the outside of the zone as injections. The PTDF represent the change in these outbound flows induced by control actions.

\subsubsection{Under normal conditions} The power flows are noted $F^N_{ij}$, with $ij \in \mathcal{Z}^L$, where  $\mathcal{Z}^L$ is the set of lines in the zone. The outbound flows $\tilde{F}^N_{kl}$ with $kl \in \mathcal{Z}^O$, the set of outbound lines, are considered as injections. $\theta_i$ with $i \in \mathcal{Z}^B$ is the phase angle and $x_{ij}$, with $x_{ij} > 0, ij \in \mathcal{Z}^L$ is the reactance of line $ij$. For the simplicity of notations, we define a battery variable $B_i$ on each node $i$ so that $B_i=B_i^{min}=B_i^{max} = 0$ if $i$ is not the battery node. $I_i$ is the net injection at node $i \in \mathcal{Z}^B$.

$\tilde{F}^{0,N}_{i}$ is the flow if no control action is taken. It is called the reference outbound flow throughout the remainder of the paper. It can be computed by a load flow on the whole grid. The flows under normal conditions are subject to the permanent ratings in \eqref{permRa}.

\begin{equation}
	\sum_{ji \in \mathcal{Z}^{L}} F_{ji}^N + \sum_{ki \in \mathcal{Z}^{O}}\tilde{F}_{ki}^{N}= I_i-B_i-C_i, \quad \forall i \in \mathcal{Z}^{B}
	\label{bilan}
\end{equation}
\begin{equation}
	F^N_{ij} = \frac{1}{x_{ij}}\cdot(\theta_i^N - \theta_j^N), \quad \forall ij \in \mathcal{Z}^{L}
	\label{phase}
\end{equation}
\begin{equation}
	\tilde{F}^N_{ij} = \tilde{F}^{0,N}_{ij} - \sum_{k \in \mathcal{Z}^{B}} ptdf(k,ij) \cdot (B_k + C_k), \quad \forall ij \in \mathcal{Z}^{O}
	\label{ptdf}
\end{equation}

\begin{equation}
	|F^N_{ij}| \leq  \bar{F_{ij}}^{perm}, \quad \forall ij \in \mathcal{Z}^{L}
	\label{permRa}
\end{equation}

\eqref{bilan} is the classical balance equation for DC power flow. \eqref{phase} is the classical DC flow definition and \eqref{ptdf} the classical DC flow definition with the PTDF formulation.

\subsubsection{Under a contingency - immediate rating}

Power flows $F^C_{ij,n}$ and phase angle variables $\theta^C_{i,n}$ are considered for each contingency $n \in \mathcal{Z}^{C}$. The equations are similar to the previous ones. The difference relies in the $ptdf$ coefficients (which depend on the topology), in the reference outbound flows $\tilde F_i^{0,C}$, and in the $x_{ij}$ if the contingency happens in the zone. The rating considered in this case is the immediate rating.

\begin{equation}
	\sum_{ji \in \mathcal{Z}^{L}} F_{ji,n}^{C} + \sum_{ki \in \mathcal{Z}^{O}}\tilde{F}_{ki,n}^{C}= I_i-B_i-C_i, \quad \forall i \in \mathcal{Z}^{B}, \forall n \in \mathcal{Z}^{C}
\end{equation}

\begin{equation}
	F^{C}_{ij,n} = \frac{1}{x_{ij}}\cdot(\theta_{i,n}^{C} - \theta_{j,n}^{C}), \quad \forall ij \in \mathcal{Z}^{L}, \forall n \in \mathcal{Z}^{C}
\end{equation}
\begin{multline}
	\tilde{F}^{C}_{ij,n} = \tilde{F}^{0,C}_{ij,n} - \sum_{k \in \mathcal{Z}^{B}} ptdf(k,ij) \cdot (B_k + C_k), \\\forall ij \in \mathcal{Z}^{O}, \forall n \in \mathcal{Z}^{C}
\end{multline}
\begin{equation}
	|F_{ij,n}^{C}| \leq \bar{F_{ij}}^{imm}, \quad \forall ij \in \mathcal{Z}^{L}, \forall n \in \mathcal{Z}^{C}
\end{equation}

\subsubsection{Under a contingency - long-term rating}
These constraints must ensure that, in the case of a contingency, the fast curative actions are able to reduce the flows under the long-term rating. The difference between these constraints and the previous ones is the presence of the curative battery injections $B^{cur}_{i,n}$ and the lower ratings $\displaystyle \bar{F_{ij}}^{long}$. The power flow and phase angle variables after the use of the fast curative actions are noted $F^{'C}_{ij,n}$ with $ij \in \mathcal{Z}^{L}, n \in \mathcal{Z}^{C}$ and $\theta^{'C}_{i,n}$ with $i \in \mathcal{B}^{L}, n \in \mathcal{Z}^{C}$.

{\small
\begin{multline}
	\sum_{ji \in \mathcal{Z}^{L}} F_{ji,n}^{'C} + \sum_{ki \in \mathcal{Z}^{O}} \tilde{F}_{ki,n}^{'C}= I_i -(B_i+B_{i,n}^{cur})-C_i, \\ \forall i \in \mathcal{Z}^{B}, \forall n \in \mathcal{Z}^{C}
\end{multline}
}
{\small
\begin{equation}
	F^{'C}_{ij,n} = \frac{1}{x_{ij}}\cdot(\theta_{i,n}^{'C} - \theta_{j,n}^{'C}), \quad \forall ij \in \mathcal{Z}^{L},  \forall n \in \mathcal{Z}^{C}
\end{equation}
}
{\small
\begin{multline}
	\tilde{F}^{C}_{ij,n} = \tilde{F}^{0,C}_{ij,n} -  \sum_{k \in \mathcal{Z}^{B}} ptdf(k,ij) \cdot (B_k + B_{k,n}^{cur}+ C_k), \\ \forall ij \in \mathcal{Z}^{O}, \forall n \in \mathcal{Z}^{C}
\end{multline}
}
{\small
\begin{equation}
	|F_{ij,n}^{'C}| \leq \bar{F_{ij}}^{long}, \quad \forall ij \in \mathcal{Z}^{L}, \forall n \in \mathcal{Z}^{C}
\end{equation}
}

\subsubsection{Under a contingency - permanent rating}
These constraints must ensure that, in case of a contingency, the curative actions (fast and long) are able to reduce the flows under the long-term rating. The difference between these constraints and the previous ones is the presence of the curative renewable generation $C^{cur}_{i,n}$ and the lower rating $ \displaystyle \bar{F_{ij}}^{perm}$. The power flow and phase angle variables after the use of all curative actions are noted $F^{''C}_{ij,n}$ with $ij \in \mathcal{Z}^{L}, n \in \mathcal{Z}^{C}$ and $\theta^{''C}_{i,n}$ with $i \in \mathcal{B}^{L}, n \in \mathcal{Z}^{C}$.

{\small \begin{multline}
	\sum_{ji \in \mathcal{Z}^{L}} F_{ji,n}^{''C} + \sum_{ki \in \mathcal{Z}^{O}} \tilde{F}_{ki,n}^{''C}= I_i-(B_i+B_{i,n}^{cur})- (C_i+C_{i,n}^{cur}), \\
	 \forall i \in \mathcal{Z}^{B}, \forall n \in \mathcal{Z}^{C}
\end{multline}}
{\small \begin{equation}
	F^{''C}_{ij,n} = \frac{1}{x_{ij}}\cdot(\theta_{i,n}^{''C} - \theta_{j,n}^{''C}), \quad \forall ij \in \mathcal{Z}^{L}, \forall n \in \mathcal{Z}^{C}
\end{equation}}
{\small \begin{multline}
	\tilde{F}^{''C}_{ij,n} = \tilde{F}^{''0,C}_{ij,n} - \sum_{k \in \mathcal{Z}^{B}} ptdf(k,ij) \cdot (B_k + B_{k,n}^{cur}+ C_k + C_{k,n}^{cur}), \\\forall ij \in \mathcal{Z}^{O}, \forall n \in \mathcal{Z}^{C}
\end{multline}}
{\small \begin{equation}
	|F_{ij,n}^{''C}| \leq \bar{F_{ij}}^{perm}, \quad \forall ij \in \mathcal{Z}^{L}, \forall n \in \mathcal{Z}^{C}
\end{equation}}

\subsubsection{Bounds on the control actions} The battery injections and the renewable curtailment have to be within their limits:
\begin{equation}
	B_i^{min} \leq B_i \leq  B_i^{max}
\end{equation}
\begin{equation}
	B_i^{min} \leq B_i + B_{i,n}^{cur}\leq  B_i^{max}, \quad \forall n \in \mathcal{Z}^{C}
\end{equation}
\begin{equation}
	0 \leq C_i \leq  C_i^{max}, \quad \forall i \in \mathcal{Z}^{B}
\end{equation}
\begin{equation}
	0 \leq C_i + C^{cur}_{i,n} \leq  C_i^{max}, \quad \forall i \in \mathcal{Z}^{B}, \forall n \in \mathcal{Z}^{C}
	\label{bounds}
\end{equation}

The bandwidths in power are obtained by the solution of the two problems with the objective function \eqref{objective2} and the constraints \eqref{bilan}-\eqref{bounds}. Indeed, the lower and upper bounds are given by the value of the preventive injection $B$ in the optimal solution. They are noted  $\underline{B}^t$ and $\displaystyle \bar{B}^t$, with $t \in \{0,T\}$, the horizon considered. The values of $B^{cur}_{i,n}$, correspond to the corrective actions and are not included in the bandwidths in power.  However, $B^{cur}_{i,n}$ is taken into account during the computation of bandwidths in energy.

\subsection{Computation of the bandwidths in energy}

The bandwidths in energy take into account the time-dependence of the problem and, by anticipating congestions, determine conditions on the battery state of charge to be able to respect the bandwidths in power. More precisely, if, at the beginning of the first time step, the battery state of charge is within the lower and upper bounds of the bandwidths in energy, then it exists a battery trajectory such that the power injections are within the bandwidths in power and the state of charge is within the bandwidths in energy for the following time steps. The bandwidths in energy represent the margins that have to be taken from the battery maximum and minimum states of charge to be able to charge or discharge according to the bandwidths in power. They can be obtained by backward computations from the bandwidths in power. The upper bound of the bandwidth is given by:

{\small \begin{equation*}
	\bar{SC}^T = SC^{max} 
\end{equation*}}
{\small \begin{equation*}
	\bar{SC}^{t} =  \bar{SC}^{t+1} - \Delta t \cdot \underline{B}^t - \Delta t_{cur} \cdot \max_{n \in  \mathcal{Z}^C} B^{cur,t}_{n}, \quad \forall t \in \{0,T-1\}
\end{equation*}}

$\Delta t$ corresponds to the time step of the forecasted data. For instance, RTE has forecasted data for each hour, so $\Delta t~=~1h$. $\Delta t_{cur}\cdot\max_{n \in \mathcal{Z}^C}~B^{cur,t}_n$ is the energy necessary to relieve the worst contingency congestion. $\Delta t_{cur}$ is the time necessary to curtail renewable generation. For RTE, renewable generation curtailment can be implemented usually within an approximate 5 minutes duration. The amount of energy required in case of a contingency to decrease flows under the long-term rating is consequently not very important. The lower bounds of the bandwidths in energy are computed in a similar way.

\section{Simulations and numerical results}

The framework for battery scheduling allowing coordination between congestion management and other services has been simulated on one of the three zones of the RTE sub-transmission grid selected for the installation of a battery in 2020. The bandwidths have been computed for a whole year using RTE hourly data with projected installed wind capacity in 2020. The time step considered in the simulations is thus one hour.

\subsection{A 90kV zone}

\begin{figure}[!h]
	\centering
	\includegraphics[scale=0.55]{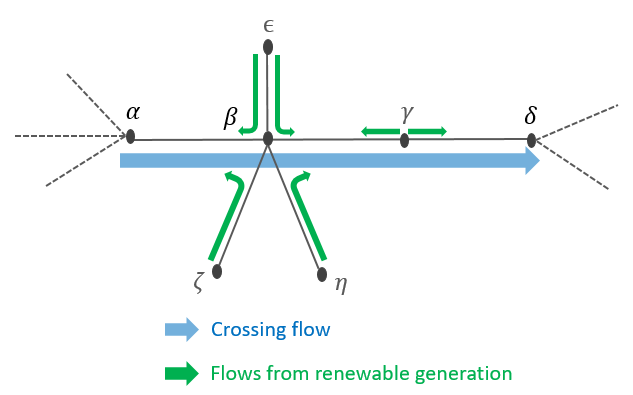}
	\caption{Map of the 90 kV zone}
	\label{map}
\end{figure}

Fig. \ref{map} shows the 90 kV studied zone, located in the West of France. The battery (12MW-24MWh) is located at substation~$\gamma$. High wind generation capacities are located on the left of $\alpha$, creating an important crossing flow between $\alpha$ and $\delta$. Table \ref{statsProd} presents the average and maximum crossing flows between $\alpha$ and $\delta$, as well as the load and generation in the zone. Renewable generation is present on each node of the zone and can be curtailed to reduce the flows below the ratings if the battery capacity is not sufficient. Table \ref{ratings} contains the lines ratings.

\begin{table}[!h]
\renewcommand{\arraystretch}{1.3}
\label{statsProd}
\centering
\begin{tabular}{|c||c|c|}
\cline{2-3} \multicolumn{1}{c|}{} & Average & Maximum\\
\hline
Generation & 49.4 & 208.0\\
\hline
Load & 10.1 & 22.3\\
\hline
Crossing flows & 34.6 & 83.2\\
\hline
\end{tabular}
\vspace{2mm}
\caption{Load, generation and crossing flows (MW)}
\end{table}

\begin{table}[!h]
\renewcommand{\arraystretch}{1.3}
\label{ratings}
\centering
\begin{tabular}{|c||c|c|c|c|}
\cline{2-4} \multicolumn{1}{c|}{}  & \multicolumn{3}{|c|}{$\alpha-\beta$ and $\beta-\gamma$} \\
\cline{2-4} \multicolumn{1}{c|}{}  & Perm. rating & Long-term rating & Imm. rating\\
\hline
Summer & 70 & 81 & 101\\
\hline
Winter& 81 & 99 & 101\\
\hline
\multicolumn{4}{c}{} \\
\cline{2-4} \multicolumn{1}{c|}{}  & \multicolumn{3}{|c|}{$\gamma-\delta$} \\

\cline{2-4} \multicolumn{1}{c|}{} & Perm. rating & Long-term rating & Imm. rating\\
\hline
Summer & 77 & 82 & 111\\
\hline
Winter& 87 & 100 & 111\\
\hline
\end{tabular}
\vspace{2mm}
\caption{Seasonal ratings in MW}
\end{table}

\subsection{Congestion under normal conditions}
Forecasted flows may exceed permanent ratings creating congestions even without any contingency. Control actions must be taken preventively to avoid exceeding the limit during operation. In the studied zone presented in Fig. \ref{map}, congestion under normal conditions can occur when the crossing flow and the renewable generation between $\alpha$ and $\delta$ are high. This type of congestions is then mostly situated on the $\gamma-\delta$ line, which has to convey the most important flow (see Fig. \ref{map}).

\begin{figure}[!h]
	\centering
	\includegraphics[scale=0.32]{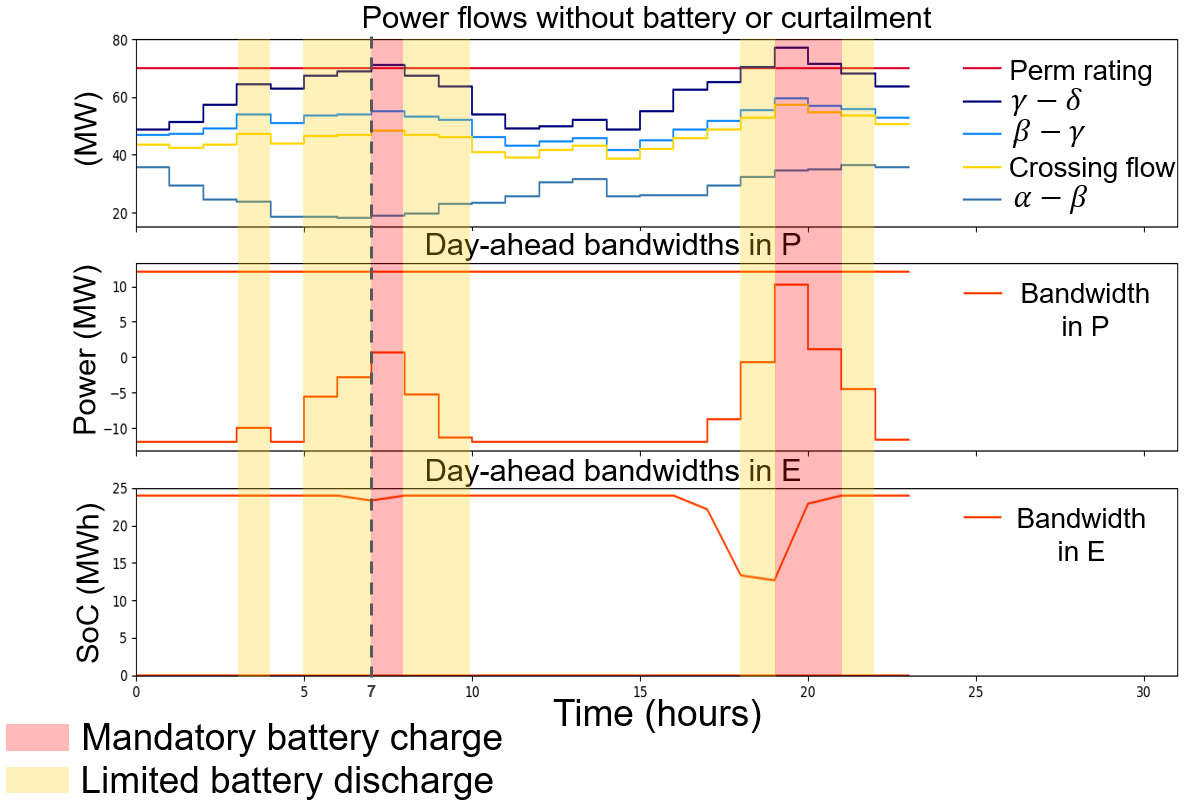}
	\caption{Bandwidths reduced because of congestions under normal conditions}
	\label{congestionN}
\end{figure}

Fig. \ref{congestionN} shows the bandwidths obtained during a summer day where there are congestions under normal conditions. The red zones indicate time steps  during which the $\gamma-\delta$ flow would be over the permanent rating without control (see first graph). Thus, the battery has to charge during these times steps to relieve the congestion. The lower bounds of the bandwidth in power are strictly positive (see second graph). The bandwidths in energy anticipate the battery mandatory charge by reducing the upper bounds on the state of charge at the beginning of each time step where the battery has to charge (see third graph). For instance, on time step 7, the permanent rating on line $\gamma-\delta$ is exceeded by 1MW. As the PTDF of line $\gamma$ on $\gamma-\delta$ is 0.6, the battery has to charge at least 1.67MW to reduce the flow by 1MW on line $\gamma-\delta$. Consequently, the bandwidth in power is [1.67MW,12MW]. The interval in which the battery state of charge must be at the beginning of the time step is [0,23.33].

The yellow zones indicate hours during which the battery discharge is limited to ensure that the battery operation will not create a congestion.

\subsection{Congestion under contingency}
In the case of a contingency, a congestion can arise under two conditions:
\begin{itemize}
	\item The forecasted flows are expected to exceed the immediate rating if a contingency occurs.
	\item The forecasted flows are expected to remain below the immediate rating, but will exceed the permanent rating in case of a contingency.
\end{itemize}

\begin{figure}[!h]
	\centering
	\includegraphics[scale=0.3]{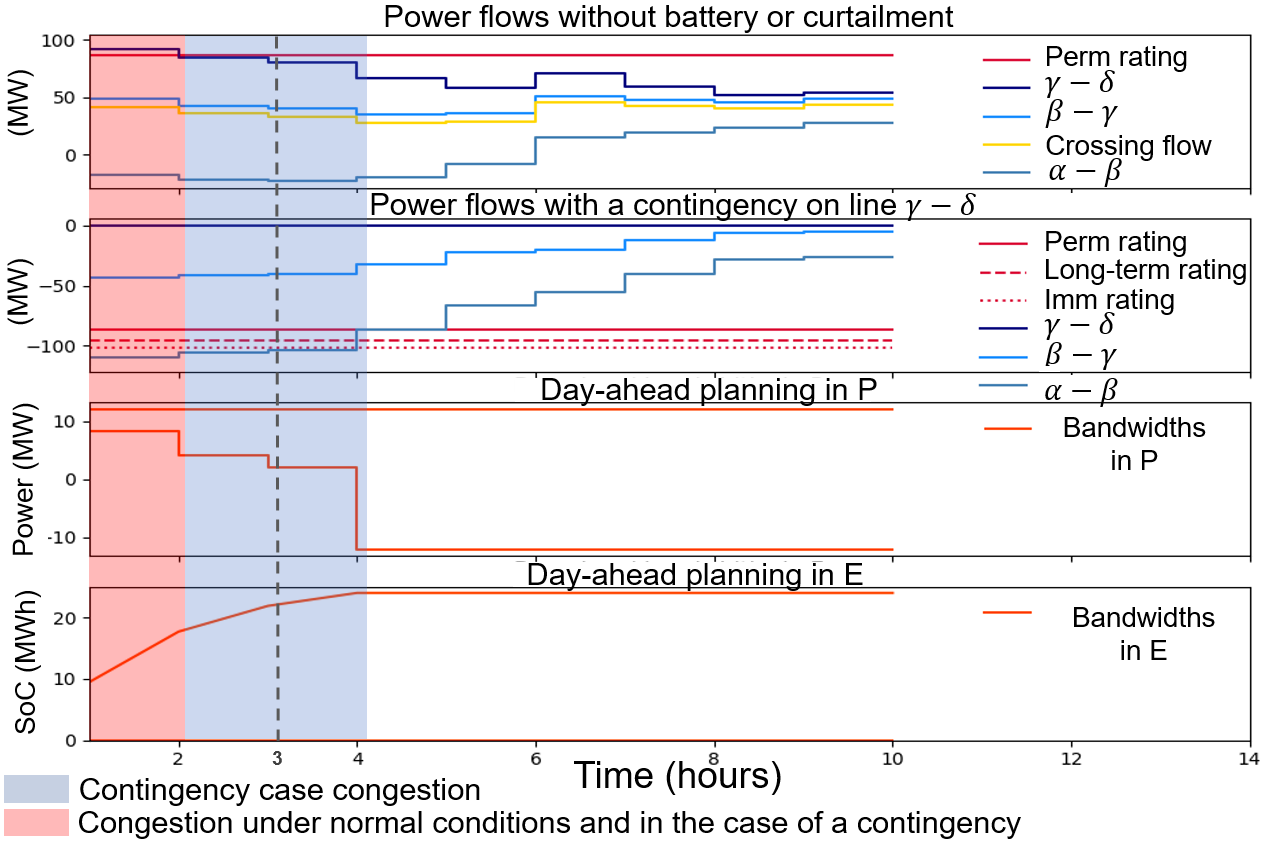}
	\caption{Bandwidths reduced because of congestions under contingency}
	\label{congestionN_1}
\end{figure}

The first case is illustrated in Fig. \ref{congestionN_1}. The blue zone indicates that there is a congestion in case of contingency. The first graph shows that under normal conditions, permanent ratings are not exceeded. The second graph shows that the $\alpha-\beta$ immediate rating would be exceeded if a contingency occurs on $\gamma-\delta$ line. Thus, the battery has to charge preventively to avoid this situation: the bandwidth in power is reduced. The bandwidth in energy is reduced accordingly. For instance on time step 3, if a contingency occurs on line $\gamma-\delta$, the flow on $\alpha-\beta$ would exceed the immediate rating by 3MW. As the PTDF of $\gamma$ on the line $\alpha-\beta$ is 1 because of a contingency on the line $\gamma-\delta$, the battery has to charge at least 3MW during the time step. Thus, the battery needs to charge at least 3MW during the time step. The bandwidth in power is [3MW,12MW], and the battery state of charge must belong to the interval [0MWh, 21MWh] at the beginning of the time step 3.

The red zone indicates two time steps during which there are both a congestion under normal conditions and in the case of a contingency. Under normal conditions, the power flow on $\gamma-\delta$ is 81 MW, which exceeds the permanent rating by 4MW. The PTDF of $\gamma$ node on line $\gamma-\delta$ is 0.6 under normal conditions. Thus, the battery has to charge 6.67MW to relieve this congestion. In the contingency case, the power flow on $\alpha-\beta$ is -110MW, which is higher in absolute value by 9MW compared to the immediate rating (101MW). The battery has to charge 9MW (PTDF=1 in the contingency case) to relieve the congestion. The bandwidth in power is thus [9MW,12MW] during time steps 0 and 1.

\subsection{Results}

Simulations conducted on the whole year give an idea of what the battery utilization rate for congestion management may be. Table \ref{year} summarizes the battery availability according to seasons. The terminology "strong congestion" refers to the case in which the battery has to charge at its maximal capacity, and is consequently not available for other services. More generally, "congestion" includes all strong congestions, but also reduced bandwidths. When the congestion is not critical, some capacities can still be shared with other services.

\begin{table}[!h]
\renewcommand{\arraystretch}{1.3}
\label{year}
\centering
\begin{tabular}{|c||c|c|c|}
\cline{2-4} \multicolumn{1}{c|}{} & Strong congestion & Congestion & \textbf{Batt. fully available}\\
\hline
Summer & $3\%$ & $6\%$ & $\bm{94\%}$\\
\hline
Winter& $9\%$& $15\%$ & $\bm{85\%}$\\
\hline
\end{tabular}
\vspace{2mm}
\caption{Congestion occurences}
\end{table}

As expected, the battery is more used in winter. Indeed, even if the ratings are higher, the wind is usually more important in this area during winter, \textit{i-e} the flows are higher.
The results show that the battery is often available for other services and an additional value can be created under multi-service.

\section{Conclusion}
During the ongoing three years of RTE's demonstration, the amount of energy stored in the three batteries will be constant \cite{decisionInvestissementRINGO}, in order to interfere with the electricity market in the same way than power lines. This paper presents the scheduling framework that will be used by RTE for the day-ahead battery operation planning. The three batteries will manage local congestions and their injections will be balanced within the bandwidths. This guarantees that the balancing constraint will not create congestions in other zones. Furthermore, the framework can be used to coordinate the battery operation for congestion management with other services. After the three-year test period, the bandwidths might be used by other services, such as frequency regulation or energy arbitrage. 

The framework presented in this paper does not consider uncertainty. Indeed, the bandwidths are computed under a deterministic set up. However, uncertainty may affect the residual capacity which may be an issue depending on the firmness of the other services provided by the battery. The impact of uncertainty on the bandwidths calculation will be studied in the future. Moreover, the batteries we consider in this paper are at the sub-transmission level and help managing local congestions. If more batteries are connected to the grid, they may be used for congestion management on the higher voltage level grid and thus help both local and interregional congestion management. At that time the global scheduling process will be built upon the work presented here and the work in \cite{wen2015enhanced}.






\bibliographystyle{IEEEtran}
\bibliography{mybib}

\end{document}